\input Tex-document.sty

\input psfig.sty

\def\volumeno{I}
\def\firstpage{159}
\def\lastpage{178}
\pageno=159
\def\shorttitle{Hyperbolic Systems of Conservation Laws}
\def\shortauthor{A. Bressan}

\title{\centerline{Hyperbolic Systems of Conservation Laws}
\centerline{in One Space Dimension}}

\author{Alberto Bressan\footnote{\eightrm *}{\eightrm S.I.S.S.A., Via
Beirut 4, Trieste 34014, Italy. E-mail: bressan@sissa.it}}

\vskip 7mm

\centerline{\boldnormal Abstract}

\vskip 4.5mm

{\narrower \ninepoint \smallskip Aim of this paper is to review
some basic ideas and recent developments in the theory of strictly
hyperbolic systems of conservation laws in one space dimension.
The main focus will be on the uniqueness and stability of entropy
weak solutions and on the convergence of vanishing viscosity
approximations.

\vskip 4.5mm

\smallskip


}

\vskip 9mm

\head{1. Introduction}

By a
{\it system of conservation laws} in $m$ space dimensions
we mean a first order system of
partial differential equations in divergence form:
$${\partial\over\partial t} U+\sum_{\alpha=1}^m
{\partial\over\partial x_\alpha}
F_\alpha(U)=0\,,\qquad\qquad U\in{I\!\!R}^n,~~(t,x)\in{I\!\!R}\times{I\!\!R}^m.$$
The components of the vector $U=(U_1,\ldots,U_n)$ are the
{\it conserved quantities}.
Systems of this type express the balance equations of continuum physics,
when small dissipation effects are neglected.
A basic example
is provided by the equations of non-viscous gases,
accounting for the conservation of mass, momentum and energy.
The subject is thus very classical, having a long tradition
which can be traced back to Euler (1755) and includes
contributions by Stokes, Riemann, Weyl and Von Neumann, among
several others.
The continued attention of analysts and mathematical physicists
during the span of over two centuries, however, has not accounted
for a comprehensive mathematical theory. On the contrary, as remarked
in [Lx2], [D2], [S2], the field is still replenished with
challenging open problems.
In several space dimensions, not even the global existence of solutions
is presently known, in any significant degree of generality.
Until now, most of the analysis has been concerned with the
one-dimensional case, and it is only here that
basic questions could be settled.
In the remainder of this paper we shall thus consider systems
in one space dimension, referring to
the books of Majda [M], Serre [S1] or Dafermos [D3] for a discussion of
the multidimensional case.

Toward a rigorous mathematical analysis of solutions, the main
difficulty that one encounters is the lack of regularity. Due to the
strong nonlinearity of the equations and the absence of diffusion terms
with smoothing effect, solutions which are initially smooth may become
discontinuous within finite time. In the presence of discontinuities,
most of the classical tools of differential calculus do not apply.
Moreover, for general $n\times n$ systems, the powerful techniques of
functional analysis cannot be used. In particular, solutions cannot be
represented as fixed points of a nonlinear transformation, or in
variational form as critical points of a suitable functional.  Dealing
with vector valued functions, comparison arguments based on upper and
lower solutions do not apply either. Up to now, the theory of
conservation laws has progressed largely by {\it ad hoc} methods.  A
survey of these techniques is the object of the present paper. \vskip
0.8em The Cauchy problem for a system of conservation laws in one space
dimension takes the form
$$u_t+f(u)_x=0,\eqno(1.1)$$
$$u(0,x)=\bar u(x).\eqno(1.2)$$
Here $u=(u_1,\ldots,u_n)$ is the vector of {\it conserved quantities},
while the components of $f=(f_1,\ldots,f_n)$ are the {\it fluxes}.
We shall always
assume that the flux function $f:{I\!\!R}^n\mapsto{I\!\!R}^n$ is smooth and that
the system is {\it
strictly hyperbolic}, i.~e., at each point $u$
the Jacobian matrix $A(u)=Df(u)$ has $n$ real, distinct eigenvalues
$$\lambda_1(u)<\cdots <\lambda_n(u).\eqno(1.3)$$
As already mentioned, a
distinguished feature of nonlinear hyperbolic systems is
the possible loss of regularity.  Even with smooth initial data,
it is well known that the solution can develop shocks in finite time.
Therefore, solutions defined globally in time can only be
found within a space of discontinuous functions.
The equation (1.1) must then be interpreted in distributional sense.
A vector valued function $u=u(t,x)$ is a {\it weak solution} of
(1.1) if
$${\int\!\!\int} \big[u\,\phi_t+f(u) \,\phi_x\big]\,dxdt=0\eqno(1.4)$$
for every test function $\phi\in{\cal C}^1_c$, continuously differentiable with
compact support.
In particular, the piecewise constant function
$$u(t,x)\doteq\cases{u^-\quad &if\quad $x<\lambda t\,$,\cr
u^+\quad &if\quad $x>\lambda t\,$,\cr}\eqno(1.5)$$
is a weak solution of (1.1) if and only if the left and right states
$u^-,u^+$ and the speed $\lambda$ satisfy the
famous Rankine-Hugoniot equations
$$f(u^+)-f(u^-)=\lambda\,(u^+-u^-)\,.\eqno(1.6)$$

When discontinuities are present, the weak solution of a Cauchy problem
may not be unique. To single out a unique ``good'' solution, additional
{\it entropy conditions} are usually imposed along shocks [Lx1], [L3].
These conditions often have a physical motivation, characterizing those
solutions which can be recovered from higher order models, letting the
diffusion or dispersion coefficients approach zero (see [D3]). \vskip
0.8em In one space dimension, the mathematical theory of hyperbolic
systems of conservation laws has developed along two main lines. \vskip
0.8em \noindent{\bf 1.} The $BV$ setting, pioneered by Glimm (1965).
Solutions are here constructed within a space of functions with bounded
variation, controlling the $BV$ norm by a wave interaction potential.
\vskip 0.8em \noindent{\bf 2.} The ${\bf L}^\infty$ setting, introduced
by DiPerna (1983), based on weak convergence and a compensated
compactness argument. \vskip 0.8em Both approaches yield results on the
global existence of weak solutions. However, it is only in the $BV$
setting that the well posedness of the Cauchy problem could recently be
proved, as well as the stability and convergence of vanishing viscosity
approximations. On the other hand, a counterexample in [BS] indicates
that similar results cannot be expected, in general, for solutions in
${\bf L}^\infty$. In the remainder of this paper we thus concentrate on
the theory of $BV$ solutions, referring to [DP2] or [S1] for the
alternative approach based on compensated compactness.

We shall first review the main ideas involved in the
construction of weak solutions, based on the Riemann problem and
the wave interaction functional.
We then present more recent results on stability, uniqueness and
characterization of entropy weak solutions. All this material
can be found in the monograph [B3].
The last section
contains an outline of the latest work
on stability and convergence of vanishing viscosity approximations.

\head{2. Existence of weak solutions}

Toward the construction of more general solutions of (1.1), the
basic building block is the {\it Riemann problem}, i.e.~the initial
value problem where the data are piecewise
constant, with a single jump at the origin:
$$u(0,x)=\cases{u^-\quad &if\quad $x<0\,$,\cr
u^+\quad &if\quad $x>0\,$.\cr}\eqno(2.1)$$ Assuming that the amplitude
$|u^+-u^-|$ of the jump is small, this problem was solved in a classical
paper of Lax [Lx1], under the additional hypothesis \vskip 0.8em
\item{(H)} For each $i=1,\ldots,n$, the $i$-th field
is either  {\it genuinely nonlinear}, so that $D\lambda_i(u)\cdot
r_i(u)>0$ for all $u$, or {\it linearly degenerate}, with
$D\lambda_i(u)\cdot r_i(u)= 0$ for all $u$. \vskip 0.8em \noindent The
solution is self-similar: $u(t,x)=U(x/t)$. It consists of $n+1$ constant
states $\omega_0=u^-$, $\omega_1,\ldots,\omega_n=u^+$ (see Fig.~1). Each
couple of adiacent states $\omega_{i-1}$, $\omega_i$ is separated either
by a {\it shock} (the thick lines in Fig.~1) satisfying the Rankine
Hugoniot equations, or else by a {\it centered rarefaction}.  In this
second case, the solution $u$ varies continuously between $\omega_{i-1}$
and $\omega_i$ in a sector of the $t$-$x$-plane (the shaded region in
Fig.~1) where the gradient $u_x$ coincides with an $i$-eigenvector of
the matrix $A(u)$.

\vskip 8pt
\centerline{\hbox{\psfig{figure=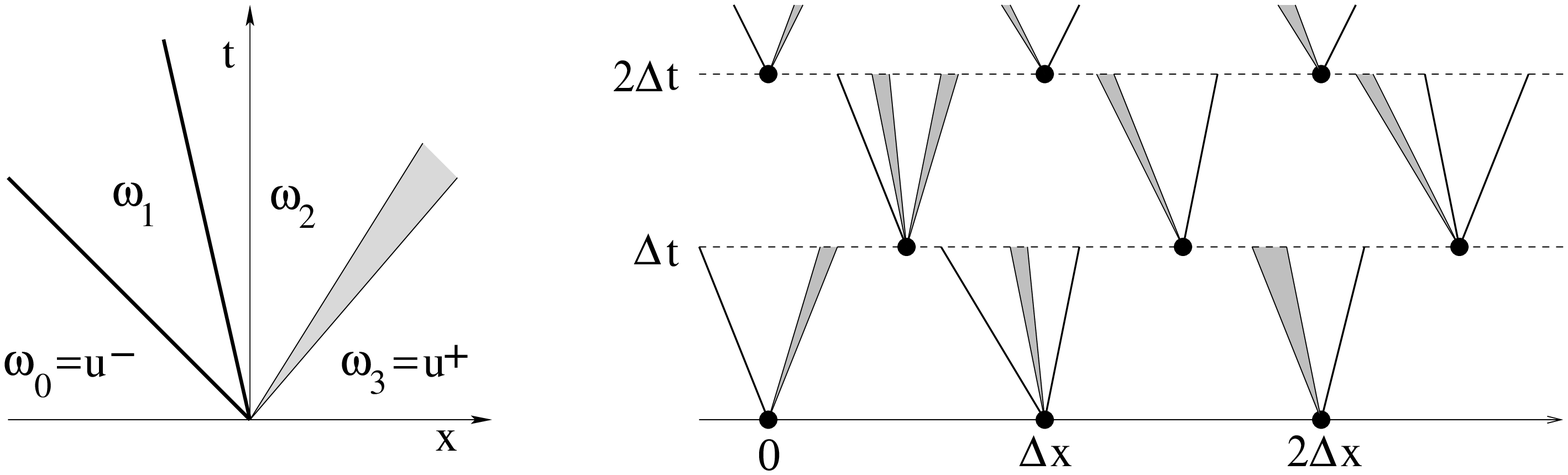,width=12cm}}}
\centerline{Figure 1~~~~~~~~~~~~~~~~~~~~~~~~ ~~~~~~~~~~~~~~~~~~~~~Figure
2~~~~~~} \vskip 8pt

Approximate solutions to a more general Cauchy problem can be
constructed by patching together several solutions of Riemann
problems. In the Glimm scheme (Fig.~2), one works with a fixed
grid in the $x$-$t$ plane, with mesh sizes $\Delta x$, $\Delta t$.
At time $t=0$ the initial data is approximated by a piecewise
constant function, with jumps at grid points. Solving the
corresponding Riemann problems, a solution is constructed up to a
time $\Delta t$ sufficiently small so that waves generated by
different Riemann problems do not interact. By a random sampling
procedure, the solution $u(\Delta t,\cdot)$ is then approximated
by a piecewise constant function having jumps only at grid points.
Solving the new Riemann problems at every one of these points, one
can prolong the solution to the next time interval $[\Delta
t,~2\Delta t]$, etc$\ldots$

\vskip 8pt
\centerline{\hbox{\psfig{figure=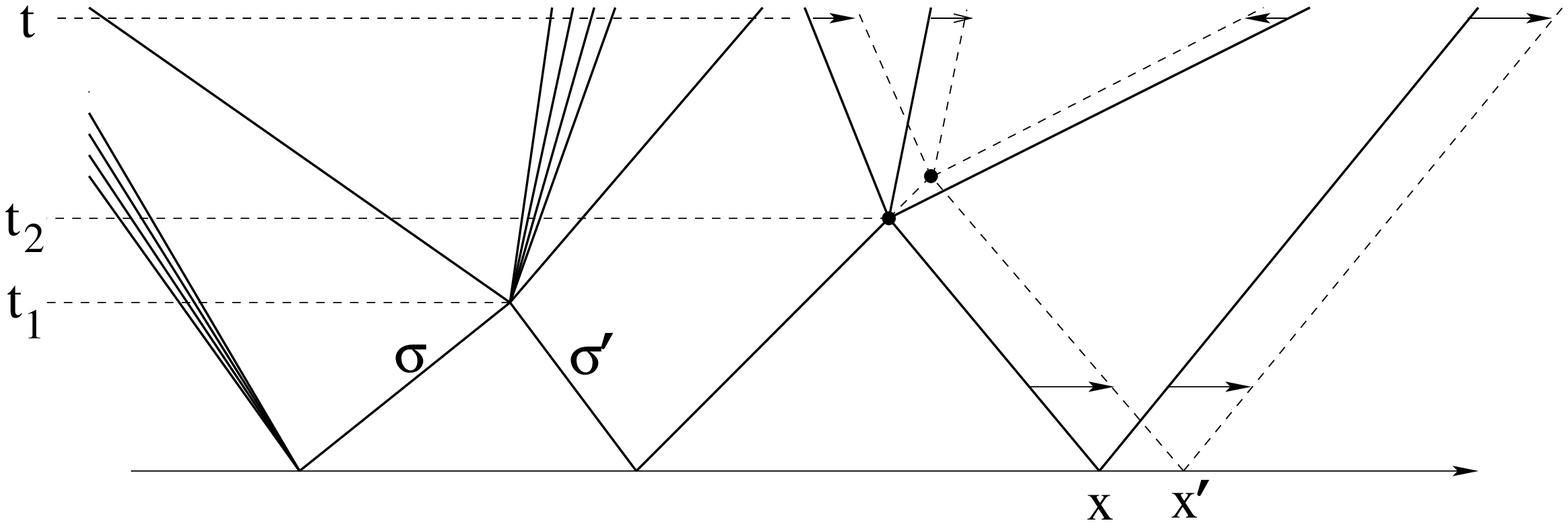,width=10cm}}}
\centerline{Figure 3} \vskip 8pt

An alternative technique for contructing approximate solutions is
by wave-front tracking (Fig.~3). This method was introduced by
Dafermos [D1] in the scalar case and later developed by various
authors [DP1], [B1], [R], [BJ]. It now provides an efficient tool
in the study of general $n\times n$ systems of conservation laws,
both for theoretical and numerical purposes [B3], [HR].

The initial data is here approximated with
a piecewise constant function, and each Riemann problem
is solved approximately, within the class of piecewise constant functions.
In particular, if the exact solution contains a centered rarefaction,
this must be approximated by a {\it rarefaction fan},
containing several small jumps.
At the first time $t_1$ where two fronts interact, the new Riemann
problem is again approximately solved by a piecewise constant function.
The solution is then prolonged
up to the second interaction time $t_2$, where the
new Riemann problem is solved, etc$\ldots$
The main difference is that in the Glimm scheme one specifies a priori
the nodal points where the the Riemann problems are to be solved.
On the other hand, in a solution constructed by wave-front tracking
the locations of the jumps and of the interaction points depend
on the solution itself, and no restarting procedure is needed.

In the end, both algorithms produce a sequence of approximate
solutions, whose convergence
relies on a compactness argument based on uniform bounds on
the total variation.
We sketch the main idea
involved in these a priori BV bounds.
Consider a piecewise constant function
$u:{I\!\!R}\mapsto {I\!\!R}^n$, say with jumps at points
$x_1<x_2<\cdots <x_N$.  Call $\sigma_\alpha$ the amplitude of the jump at
$x_\alpha$.
The {\it total strength of waves} is
then defined as
$$V(u)\doteq \sum_{\alpha} |\sigma_\alpha|. \eqno(2.2)$$
Clearly, this is an equivalent way to measure the total variation.
Along a solution $u=u(t,x)$ constructed by front tracking,
the quantity $V(t)=V\big(u(t,\cdot)\big)$
may well increase at interaction times.
To provide global a priori bounds,
following [G] one introduces
a {\it wave interaction potential}, defined as
$$Q(u)=\sum_{(\alpha,\beta)\in{\cal A}} |\sigma_\alpha\,\sigma_{\beta}|,
\eqno(2.3)$$ where the summation runs over the set ${\cal A}$ of
all couples of approaching waves. Roughly speaking, we say that
two wave-fronts located at $x_\alpha<x_\beta$ are {\it
approaching} if the one at $x_\alpha$ has a faster speed than the
one at $x_\beta$ (hence the two fronts are expected to collide at
a future time). Now consider a time $\tau$ where two incoming
wave-fronts interact, say with strengths $\sigma$, $\sigma'$ (for
example, take $\tau=t_1$ in Fig.~3). The difference between the
outgoing waves emerging from the interaction and the two incoming
waves $\sigma,\sigma'$ is of magnitude ${\cal O}(1)\cdot
|\sigma\sigma'|$. On the other hand, after time $\tau$ the two
incoming waves are no longer approaching. This  accounts for the
decrease of the functional $Q$ in (2.3) by the amount
$|\sigma\sigma'|$. Observing that the new waves generated by the
interaction could approach all other fronts, the change in the
functionals $V,Q$ across the interaction time $\tau$ is estimated
as
$$\Delta V(\tau)={\cal O}(1)\cdot|\sigma\sigma'|\,,\qquad\qquad
\Delta Q(\tau)=-|\sigma\sigma'|+{\cal O}(1)\cdot|\sigma\sigma'|\,V(\tau-).$$
If the initial data has small total variation, for a suitable
constant $C_0$ the quantity
$$\Upsilon(t)\doteq V\big(u(t,\cdot)\big)+C_0\,Q\big(u(t,\cdot)\big)$$
is monotone decreasing in time. This argument provides the uniform BV
bounds on all approximate solutions. Using Helly's compactness theorem,
one obtains the convergence of a subsequence of approximate solutions,
and hence the  global existence of a weak solution. \vskip 0.8em
\noindent{\bf Theorem 1.} {\it Let the system (1.1) be strictly
hyperbolic and satisfy the assumptions (H).  Then, for a sufficiently
small $\delta>0$ the following holds. For every initial condition $\bar
u$ with
$$\|\bar u\|_{{\bf L}^\infty}<\delta\,,\qquad\qquad \hbox{Tot.Var.}\{\bar u\}<\delta\,,
\eqno(2.4)$$ the Cauchy problem has a weak solution, defined for all
times $t\geq 0$.} \vskip 0.8em This result is based on careful analysis
of solutions of the Riemann problem and on the use of a quadratic
interaction functional (2.3) to control the creation of new waves. These
techniques also provided the basis for subsequent investigations of
Glimm and Lax [GL] and Liu [L2] on the asymptotic behavior of weak
solutions as $t\to\infty$.

\head{3. Stability}

The previous existence result relied on a compactness argument which, by
itself, does not provide informations on the uniqueness of solutions. A
first understanding of the dependence of weak solutions on the initial
data was provided by the analysis of front tracking approximations. The
idea is to perturb the initial data by shifting the position of one of
the jumps, say from $x$ to a nearby point $x'$ (see Fig.~3). By
carefully estimating the corresponding shifts in the positions of all
wave-fronts at a later time $t$, one obtains a bound on the ${\bf L}^1$
distance between the original and the perturbed approximate solution.
After much technical work, this approach yielded a proof of the
Lipschitz continuous dependence of solutions on the initial data, first
in [BC1] for  $2\times 2$ systems, then in [BCP] for general $n\times n$
systems. \vskip 0.8em \noindent{\bf Theorem 2.} {\it Let the system
(1.1) be strictly hyperbolic and satisfy the assumptions (H). Then, for
every initial data $\bar u$ satisfying (2.4) the weak solution obtained
as limit of Glimm or front tracking approximations is unique and depends
Lipschitz continuously on the initial data, in the ${\bf L}^1$
distance.} \vskip 0.8em These weak solutions can thus be written in the
form $u(t,\cdot)=S_t\bar u$, as trajecories of a semigroup $S:{\cal
D}\times[0,\infty[\,\mapsto{\cal D}$ on some domain ${\cal D}$
containing all functions with sufficiently small total variation. For
some Lipschitz constants $L,L'$ one has
$$\big\|S_t\bar u-S_s\bar v\big\|_{{\bf L}^1}\leq L\,\|\bar u-\bar v\|_{{\bf L}^1}
+L'|t-s|\,,\eqno(3.1)$$ for all $t,s\geq 0$ and initial data $\bar
u,\bar v\in{\cal D}$. \vskip 0.8em An alternative proof of Theorem 2 was
later achieved by a technique introduced by Liu and Yang in [LY] and
presented in [BLY] in its final form. The heart of the matter is to
construct a nonlinear functional, equivalent to the ${\bf L}^1$
distance, which is decreasing in time along every pair of solutions.  We
thus seek $\Phi=\Phi(u,v)$ and a constant $C$ such that
$${1\over C}\cdot
\big\|v-u\big\|_{{\bf L}^1}\,\leq \,\Phi(u,v)\,\leq\, C\cdot
\big\|v-u\big\|_{{\bf L}^1}\,,
\eqno(3.2)$$
$${d\over dt} \Phi\big(u(t),~v(t)\big)\leq 0.\eqno(3.3)$$

\vskip 8pt
\centerline{\hbox{\psfig{figure=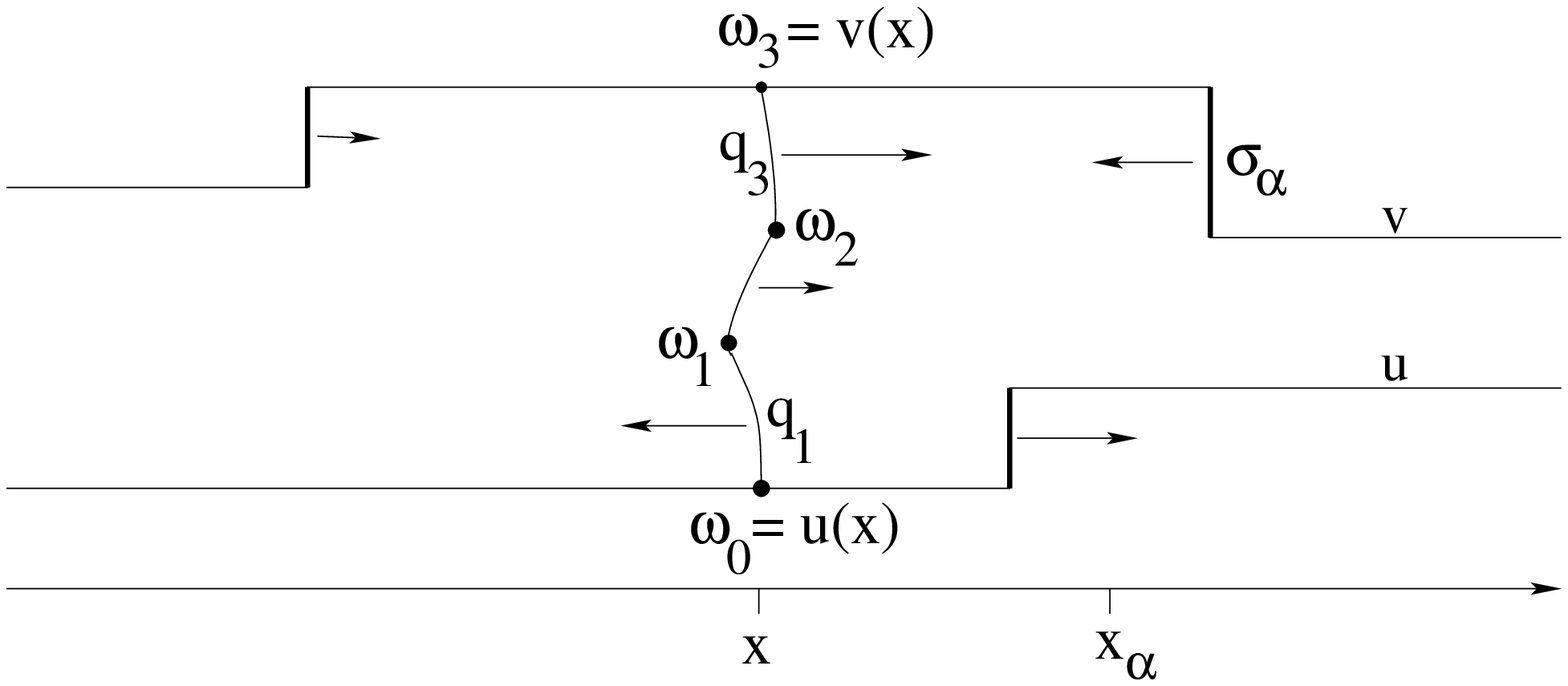,width=10cm}}}
\centerline{Figure 4} \vskip 8pt

In connection with piecewise constant functions
$u,v:{I\!\!R}\mapsto{I\!\!R}^n$ generated by a front tracking
algorithm, this functional can be defined as follows (Fig.~4). At
each point $x$, we connect the states $u(x)$, $v(x)$ by means of
$n$ shock curves.  In other words, we construct intermediate
states $\omega_0=u(x),\omega_1,\ldots,\omega_n=v(x)$ such that
each pair $\omega_{i-1},\omega_i$ is connected by an $i$-shock.
These states can be uniquely determined by the implicit function
theorem. Call $q_1,\ldots,q_n$, the strengths of these shocks. We
regard $q_i(x)$ as the $i$-th scalar component of the jump $\big(
u(x),~v(x)\big)$. For some constant $C'$, one clearly has
$${1\over C'}\cdot
\big|v(x)-u(x)\big|\leq \sum_{i=1}^n \big| q_i(x)\big|\leq  C'\cdot
\big|v(x)-u(x)\big|\,.\eqno(3.4)$$
The functional $\Phi$ is now defined as
$$\Phi(u,v)\doteq \sum_{i=1}^n \int_{-\infty}^\infty
W_i(x)\,\big|q_i(x)\big|~dx,\eqno(3.5)$$
where the weights $W_i$ take the form
$$\eqalign{ W_i(x)\doteq  1 &+ \kappa_1
\cdot\big[\hbox{total strength of waves in
$u$ and in $v$}\cr
&\qquad\qquad\qquad\qquad\hbox{which approach the $i$-wave} ~q_i(x)\big]\cr
&+\kappa_2\cdot \big[\hbox{wave interaction potentials of
$u$ and of $v$}\big]\cr
\doteq 1&+\kappa_1 V_i(x)+\kappa_2\big[ Q(u)+Q(v)\big]
\cr}\eqno(3.6)$$
for suitable constants $\kappa_1,\kappa_2$.
Notice that, by construction,
$q_i(x)$ represents the strength of a fictitious shock wave
located at $x$, travelling with a speed $\lambda_i(x)$
determined by the  Rankine-Hugoniot
equations.
In (3.6), it is thus meaningful to consider the
quantity
$$V_i(x)\doteq \sum_{\alpha\in{\cal A}_i(x)} |\sigma_\alpha|\,,$$
where the summation extends to all wave-fronts $\sigma_\alpha$ in $u$
and in $v$ which are {\it approaching} the $i$-shock $q_i(x)$. {}From
(3.4) and the boundedness of the weights $W_i$, one easily derives
(3.2). By careful estimates on the Riemann problem, one can prove that
also (3.3) is approximately satisfied. In the end, by taking a limit of
front tracking approximations, one obtains Theorem 2. \vskip 0.8em For
general $n\times n$ systems, in (3.1) one finds a Lipschitz constant
$L>1$.  Indeed, it is only in the scalar case that the semigroup is
contractive and the theory of accretive operators and abstract evolution
equations in Banach spaces can be applied, see [K], [C]. We refer to the
flow generated by a system of conservation laws as a {\it Riemann
semigroup}, because it is entirely determined by specifying how Riemann
problems are solved. As proved in [B2], if two semigroups $S,S'$ yield
the same solutions to all Riemann problems, then they coincide, up to
the choice of their domains. \vskip 0.8em {}From (3.1) one can deduce
the error bound
$$\big\| w(T)-S_T w(0)\big\|_{{\bf L}^1}\leq L\cdot\int_0^T
\left\{\liminf_{h\to 0+}{\big\| w(t +h)-S_h w(t)\big\|_{{\bf L}^1} \over
h}\right\}~dt\,,\eqno(3.7)$$ valid for every Lipschitz continuous map
$w:[0,T]\mapsto{\cal D}$ taking values inside the domain of the
semigroup. We can think of  $t\mapsto w(t)$ as an approximate solution
of (1.1), while $t\mapsto S_tw(0)$ is the exact solution having the same
initial data. According to (3.7), the distance at time $T$ is bounded by
the integral of an {\it instantaneous error rate}, amplified by the
Lipschitz constant $L$ of the semigroup. \vskip 0.8em Using (3.7), one
can estimate the distance between a front tracking approximation and the
corresponding exact solution. For approximate solutions constructed by
the Glimm scheme, a direct application of this same formula is not
possible because of the additional errors introduced by the restarting
procedures at times $t_k\doteq k\,\Delta t$. However, relying on a
careful analysis of Liu [L1], one can construct a front tracking
approximate solution having the same initial and terminal values as the
Glimm solution. By this technique, in [BM] the authors proved the
estimate
$$\lim_{\Delta x\to 0}{\big\| u^{\rm Glimm}(T,\cdot)-
u^{\rm exact}(T,\cdot)\big\|_{{\bf L}^1}\over\sqrt{\Delta x}\cdot |\ln \Delta x|}
= 0\,.\eqno(3.8)$$
In other words, letting the mesh sizes
$\Delta x,\Delta t\to 0$ while keeping their
ratio $\Delta x/\Delta t$ constant,
the ${\bf L}^1$ norm of the error in the Glimm approximate solution
tends to zero at a rate slightly slower than $\sqrt{\Delta x}$.

\head{4. Uniqueness}

The uniqueness and stability results stated in Theorem 2 refer to
a special class of weak solutions: those obtained as limits of
Glimm or front tracking approximations.  For several applications,
it is desirable to have a uniqueness theorem valid for general
weak solutions, without reference to any particular constructive
procedure. Results in this direction were proved in [BLF], [BG],
[BLe]. They are all based on the error formula (3.7). In the
proofs, one considers a weak solution $u=u(t,x)$ of the Cauchy
problem (1.1)--(1.2).  Assuming that $u$ satisfies suitable
entropy and regularity conditions, one shows that
$$
\liminf_{h\to 0+}{\big\| u(t +h)-S_h u(t)\big\|_{{\bf L}^1} \over
h}=0\eqno(4.1)$$ at almost every time $t$. By (3.7), $u$ thus coincides
with the semigroup trajectory $t\mapsto S_tu(0)=S_t\bar u$. Of course,
this implies uniqueness. As an example, we state below the result of
[BLe]. Consider the following assumptions: \vskip 0.8em
\item{\bf (A1)} {\bf (Conservation Equations)} The function $u=u(t,x)$
is a weak solution of the Cauchy problem (1.1)--(1.2), taking
values within the domain ${\cal D}$ of the semigroup $S$. More
precisely, $u:[0,T]\mapsto{\cal D}$ is continuous w.r.t.~the ${\bf
L}^1$ distance. The initial condition (1.2) holds, together with
$${\int\!\!\int} \big[u\,\phi_t+f(u) \,\phi_x\big]\,dxdt=0$$
for every ${\cal C}^1$ function $\phi$ with compact support contained
inside the open strip $\, ]0,T[\,\times{I\!\!R}$. \vskip 0.8em
\item{\bf (A2)} {\bf (Lax Entropy Condition)}
Let $u$ have an approximate jump
discontinuity at some point $(\tau,\xi)\in \, ]0,T[\times{I\!\!R}$.
In other words, assume that there exists states
$u^-,u^+\in\Omega$ and a speed $\lambda\in{I\!\!R}$ such that, calling
$$U(t,x)\doteq
\cases{ u^-\qquad &if\qquad $x<\xi+\lambda(t-\tau)$,
\cr u^+\qquad &if\qquad $x>\xi+\lambda(t-\tau)$,\cr}\eqno(4.2)$$
there holds
$$\lim_{\rho\to 0+}~{1\over\rho^2}
\int_{\tau-\rho}^{\tau+\rho}\int_{\xi-\rho}^{\xi
+\rho} \Big| u(t,x)-U(t,x)\Big|~dxdt=0.\eqno(4.3)$$
Then, for some $i\in\{1,\ldots,n\}$, one has the entropy inequality:
$$\lambda_i(u^-)\geq\lambda\geq\lambda_i(u^+).\eqno(4.4)$$
\vskip 0.8em
\item{\bf (A3)} {\bf (Bounded Variation Condition)}
The function $x\mapsto u\big( \tau(x), x)$ has bounded variation along
every Lipschitz continuous space-like curve $\big\{t=\tau(x)\big\}$,
which satisfies $|d\tau/dx|<\delta$ a.e., for some constant $\delta>0$
small enough. \vskip 0.8em \noindent{\bf Theorem 3.} {\it Let $u=u(t,x)$
be a weak solution of the Cauchy problem (1.1)--(1.2) satisfying the
assumptions (A1), (A2) and  (A3).   Then
$$u(t,\cdot)=S_t\bar u
\eqno(4.5)$$ for all $t$. In particular, the solution that satisfies the
three above conditions is unique.} \vskip 0.8em An additional
characterization of these unique solutions, based on local integral
estimates, was given in [B2]. The underlying idea is as follows. In a
forward neighborhood of a point $(\tau,\xi)$ where $u$ has a jump, the
weak solution $u$ behaves much in the same way as the solution of the
corresponding Riemann problem. On the other hand, on a region where its
total variation is small, our solution $u$ can be accurately
approximated by the solution of a linear hyperbolic system with constant
coefficients.

To state the result more precisely,
we introduce some notations.
Given a function $u=u(t,x)$ and a point $(\tau,\xi)$,
we denote by  $U^\sharp_{(u;\tau,\xi)}$ the solution
of the Riemann problem
with initial data
$$ u^-=\lim_{x\to \xi-} u(\tau,x),\qquad\qquad u^+=\lim_{x\to \xi+}
u(\tau,x).\eqno(4.6)$$
In addition, we define $U^\flat_{(u;\tau,\xi)}$ as the solution
of the linear hyperbolic Cauchy problem with constant coefficients
$$w_t+\widehat A w_x=0,\qquad\qquad w(0,x)=u(\tau,x).\eqno(4.7)$$
Here $\widehat A\doteq A\big(u(\tau,\xi)\big)$.
Observe that (4.7) is obtained from the quasilinear system
$$u_t+A(u)u_x=0\qquad\qquad (A=Df)\eqno(4.8)$$
by ``freezing'' the coefficients of the matrix $A(u)$ at the point
$(\tau,\xi)$ and choosing $u(\tau)$ as initial data. A new notion of
``good solution'' can now be introduced, by locally comparing a function
$u$ with the self-similar solution of a Riemann problem and with the
solution of a linear hyperbolic system with constant coefficients. More
precisely, we say that a function $u=u(t,x)$ is a {\bf viscosity
solution} of the system (1.1) if $t\mapsto u(t,\cdot)$ is continuous as
a map with values into ${\bf L}^1_{\rm loc}$, and moreover the following
integral estimates hold. \vskip 0.8em \noindent (i)  At every point
$(\tau,\xi)$, for every $\beta'>0$ one has
$$\lim_{h\to 0+} {1\over h}\int_{\xi-\beta' h}
^{\xi+\beta' h}\Big| u(\tau+h,~x)- U^\sharp
_{(u;\tau,\xi)}(h,~x-\xi)\Big|~dx~=~0.\eqno(4.9)$$ \vskip 0.8em
\noindent (ii)  There exist constants $C,\beta>0$ such that, for every
$\tau\geq 0$ and $a<\xi<b$, one has
$$\limsup_{h\to 0+}{1\over h}\int_{a+\beta h}
^{b-\beta h}\Big| u(\tau+h,~x)- U^\flat _{(u;\tau,\xi)}(h,x)\Big|~dx\leq
C\cdot \Big(\hbox{Tot.Var.}\big\{ u(\tau);~]a,\,b[ \
\big\}\Big)^2.\eqno(4.10)$$ \vskip 0.8em As proved in [B2], this concept
of viscosity solution completely characterizes semigroup trajectories.
\vskip 0.8em \noindent{\bf Theorem 4.} {\it Let $S:{\cal D}\times
[0,\infty[\times{\cal D}$ be a semigroup generated by the system of
conservation laws (1.1). A function $u:[0,T]\mapsto{\cal D}$ is a
viscosity solution of (1.1) if and only if $u(t)=S_t u(0)$ for all $t\in
[0,T]$.}

\head{5. Vanishing viscosity approximations}

A natural conjecture is that the entropic
solutions of
the hyperbolic system (1.1) actually coincide with the
limits of solutions to the parabolic system
$$u^\varepsilon_t+f(u^\varepsilon)_x=\varepsilon\,u^\varepsilon_{xx}\,,\eqno(5.1)$$
letting the viscosity coefficient $\varepsilon\to 0$. In view of the
previous uniqueness results, one expects that the vanishing viscosity
limit should single out the unique ``good'' solution of the Cauchy
problem, satisfying the appropriate entropy conditions. In earlier
literature, results in this direction were based on three main
techniques: \vskip 0.8em \noindent{\bf 1 - Comparison principles for
parabolic equations.} For a {\it scalar} conservation law, the
existence, uniqueness and global stability of vanishing viscosity
solutions was first established by Oleinik [O] in one space dimension.
The famous paper by Kruzhkov [K] covers the more general class of ${\bf
L}^\infty$ solutions and is also valid in several space dimensions.
\vskip 0.8em \noindent{\bf 2 - Singular perturbations.} Let $u$ be a
piecewise smooth solution of the $n\times n$ system (1.1), with finitely
many non-interacting, entropy admissible shocks. In this special case,
using a singular perturbation technique, Goodman and Xin [GX]
constructed a family of solutions $u^\varepsilon$ to (5.1), with
$u^\varepsilon\to u$ as $\varepsilon\to 0$. \vskip 0.8em \noindent{\bf 3
- Compensated compactness.} If, instead of a $BV$ bound, only a uniform
bound on the ${\bf L}^\infty$ norm of solutions of (5.1) is available,
one can still construct a weakly convergent subsequence
$u^\varepsilon\rightharpoonup u$. In general, we cannot expect that this
weak limit satisfies the nonlinear equations (1.1). However, for a class
of $2\times 2$ systems, in [DP2] DiPerna showed that this limit $u$ is
indeed a weak solution of (1.1). The proof relies on a compensated
compactness argument, based on the representation of the weak limit in
terms of Young measures, which must reduce to a Dirac mass due to the
presence of a large family of entropies. \vskip 0.8em Since the main
existence and uniqueness results for hyperbolic systems of conservation
laws are valid within the space of BV functions, it is natural to seek
uniform BV bounds also for the viscous approximations $u^\varepsilon$ in
(5.1). This is indeed the main goal accomplished in [BB]. As soon as
these BV bounds are established, the existence of a vanishing viscosity
limit follows by a standard compactness argument. The uniqueness of the
limit can then be deduced from the uniqueness theorem in [BG]. By
further analysis, one can also prove the continuous dependence on the
initial data for the viscous approximations $u^\varepsilon$, in the
${\bf L}^1$ norm. Remarkably, these results are valid for general
$n\times n$ strictly hyperbolic systems, not necessarily in conservation
form. \vskip 0.8em \noindent{\bf Theorem 5.} {\it Consider the Cauchy
problem for a strictly hyperbolic system with viscosity
$$u^\varepsilon_t+A(u^\varepsilon)u^\varepsilon_x=\varepsilon\,u^\varepsilon_{xx}, \qquad\qquad u^\varepsilon(0,x)
=\bar u(x)\,.\eqno(5.2)$$
Then there exist constants $C,L,L'$ and $\delta>0$ such that the
following holds.  If
$$\hbox{Tot.Var.}\{\bar u\}<\delta\,,\qquad\qquad
\big\|\bar u(x)\|_{{\bf L}^\infty}<\delta\,,
\eqno(5.3)$$
then for each $\varepsilon>0$ the Cauchy problem (5.2)
has a unique solution $u^\varepsilon$, defined for all $t\geq 0$.
Adopting a semigroup notation, this will be written as
$t\mapsto
u^\varepsilon(t,\cdot)\doteq S^\varepsilon_t\bar u$.
In addition, one has:
$${\bf BV~ bounds :}\qquad\qquad\quad
\hbox{Tot.Var.}\big\{S_t^\varepsilon \bar u\big\}\leq C\,
\hbox{Tot.Var.}\{\bar u\}\,.
\qquad\qquad\qquad\eqno(5.4)$$
$${\bf L}^1  ~{\bf stability :}\qquad\qquad
\qquad\big\|S^\varepsilon_t\bar u-S^\varepsilon_t\bar v\big\|_{{\bf L}^1}
\leq L\,
\big\|\bar u-\bar v\big\|_{{\bf L}^1}\,,
\qquad\qquad\qquad\eqno(5.5)$$
$$\big\|S^\varepsilon_t\bar u-S^\varepsilon_s\bar u\big\|_{{\bf L}^1}
\leq L'\,
\Big(|t-s|+\big|\sqrt {\varepsilon t}-\sqrt {\varepsilon s}\,\big|\Big)\,.
\eqno(5.6)
$$
\noindent {\bf Convergence.}~ As $\varepsilon\to 0+$, the
solutions $u^\varepsilon$ converge to the trajectories of a
semigroup $S$ such that
$$\big\|S_t\bar u-S_s\bar v\big\|_{{\bf L}^1}\leq L\, \|\bar u-\bar v\|_{{\bf L}^1}
+L'\,|t-s|\,.\eqno(5.7)$$
These vanishing viscosity limits
can be regarded as the unique {\bf vanishing viscosity solutions}
of the hyperbolic Cauchy problems
$$u_t+A(u)u_x=0,\qquad\qquad u(0,x)=\bar u(x)\,.\eqno(5.8)$$
\vskip 0.8em In the conservative case where $A(u)=Df(u)$ for some flux
function $f$, the vanishing viscosity solution is a weak solution of
$$u_t+f(u)_x=0,\qquad\qquad u(0,x)=\bar u(x)\,,\eqno(5.9)$$
satisfying the Liu admissibility conditions [L3].  Moreover, the
vanishing viscosity solutions are precisely the same as the viscosity
solutions defined at (4.9)--(4.10) in terms of local integral
estimates.} \vskip 0.8em The key step in the proof is to establish a
priori bounds on the total variation of solutions of
$$u_t+A(u)u_x=u_{xx}\eqno(5.10)$$
uniformly valid for all times $t\in [0,\infty[\,$. We outline here the
main ideas. \vskip 0.8em
\item{(i)} At each point $(t,x)$ we decompose the gradient along
a suitable basis of unit vectors $\tilde r_i$, say
$$u_x=\sum v_i\tilde r_i\,.\eqno(5.11)$$

\item{(ii)} We then derive an equation describing the evolution of
these gradient components
$$v_{i,t}+(\tilde\lambda_i v_i)_x-v_{i,xx}=\phi_i\,.\eqno(5.12)$$

\item{(iii)} Finally, we show that all source terms
$\phi_i=\phi_i(t,x)$ are integrable. Hence, for all $\tau>0$,
$$\big\|v_i(\tau,\cdot)\big\|_{{\bf L}^1}\leq
\big\|v_i(0,\cdot)\big\|_{{\bf L}^1}+
\int_0^\infty\int_{{I\!\!R}}\big|\phi_i(t,x)\big|\,dxdt<\infty\,.\eqno(5.13)$$
\vskip 0.8em In this connection, it seems natural to decompose the
gradient $u_x$ along the eigenvectors of the hyperbolic matrix $A(u)$.
This approach however does NOT work.   In the case where the solution
$u$ is a travelling viscous shock profile, we would obtain source terms
which are not identically zero. Hence they are certainly not integrable
over the domain $\big\{t>0\,, x\in{I\!\!R}\big\}$.

An alternative approach, proposed by S.~Bianchini, is to decompose $u_x$
as a {\it sum of gradients of
viscous travelling waves}.
By a viscous travelling $i$-wave we mean a solution of (5.10) having the form
$$w(t,x)=U(x-\sigma t)\,,\eqno(5.14)$$
where the speed $\sigma$ is close to the $i$-th eigenvalue
$\lambda_i$ of the hyperbolic matrix $A$.
Clearly, the function $U$ must provide a solution to the
second order O.D.E.
$$U''=\big(A(U)-\sigma\big) U'.\eqno(5.15)
$$
The underlying idea for the decomposition is as follows.
At each point $(t,x)$, given $(u,u_x,u_{xx})$, we seek
travelling wave profiles $U_1,\ldots,U_n$ such that
$$U_i(x)=u(x),\qquad\qquad i=1,\ldots,n\,,\eqno(5.16)$$
$$\sum_i U_i'(x)=u_x(x)\,,\qquad\qquad
\sum_i U_i''(x)=u_{xx}(x)\,.\eqno(5.17)$$ In general, the system
of algebraic equations (5.16)--(5.17) admits infinitely many
solutions. A unique solution is singled out by considering only
those travelling profiles $U_i$ that lie on a suitable {\it center
manifold} ${\cal M}_i$. We now call $\tilde r_i$ the unit vector
parallel to $U_i'$, so that $U_i'=v_i\tilde r_i$ for some scalar
$v_i$. The decomposition (5.11) is then obtained from the first
equation in (5.17).

Toward the BV estimate, the second part of the proof consists in
deriving the equation (5.12) and estimating the integrals of the source
terms $\phi_i$.  Here the main idea is that these source terms can be
regarded as generated by wave interactions. In analogy with the
hyperbolic case considered by Glimm [G], the total amount of these
interactions can be controlled by suitable Lyapunov functionals. We
describe here the main ones. \vskip 0.8em \noindent{\bf 1.} Consider
first two independent, scalar diffusion equations with strictly
different drifts:
$$\left\{ \eqalign{
z_{t}+\big[\lambda(t,x)z\big]_x-z_{xx}&=0\,, \cr
z^*_{t}+\big[\lambda^*(t,x)z^*\big]_x-z^*_{xx}&=0\,, \cr} \right.
$$
assuming that
$$\inf_{t,x}\lambda^*(t,x)-\sup_{t,x}\lambda(t,x)\geq c>0\,.$$
We regard
$z$ as the density of waves with a slow speed $\lambda$
and
$z^*$ as the density of waves with a fast speed $\lambda^*$.
A {\it transversal interaction potential} is defined as
$$Q(z,z^*) \doteq {1\over c}
\int\!\!\int_{{I\!\!R}^2} K(x_2-x_1) \big| z(x_1) \big| \,
\big|z^*(x_2)
\big|\, dx_1 dx_2\,,\eqno(5.18)$$
$$K(y) \doteq \cases{e^{-cy/2}\qquad &if\qquad $y>0\,$,\cr
~1\qquad &if\qquad $y\leq 0\,$.\cr}
\eqno(5.19)$$
One can show that this functional $Q$ is monotonically decreasing along
every couple of solutions $z,z^*$.
The total amount of interaction between fast and slow waves can now be
estimated as
$$\eqalign{\int_0^\infty\int_{{I\!\!R}}
\big|z(t,x)\big|\, \big|z^*(t,x)\big| &\,dxdt~\leq~ -\int_0^\infty
\left[{d \over dt}Q\big(z(t),\,z^*(t)\big)\right]\,dt\cr \leq
Q\big(z(0),\,z^*(0)\big)&~\leq ~ {1 \over c}
\int_{{I\!\!R}}\big|z(0,x)\big|\,dx\cdot
\int_{{I\!\!R}}\big|z^*(0,x)\big|\, dx\,. \cr}$$ By means of Lyapunov
functionals of this type one can control all source terms in (5.12) due
to the interaction of waves of different families. \vskip 0.8em
\noindent{\bf 2.} To control the interactions between waves of the same
family, we seek functionals which are decreasing along every solution of
a scalar viscous conservation law
$$u_t+g(u)_x=u_{xx}\,.\eqno(5.20)$$
For this purpose, to a scalar function $x\mapsto u(x)$ we associate the curve
in the plane
$$\gamma\doteq\pmatrix{u\cr g(u)-u_x\cr}~=~
\pmatrix{\hbox{conserved quantity}\cr\hbox{flux}\cr}.\eqno(5.21)$$
In connection with a solution $u=u(t,x)$ of (5.20),
the curve $\gamma$ evolves according to
$$\gamma_t+g'(u)\gamma_x=\gamma_{xx}\,.\eqno(5.22)$$
Notice that the vector $g'(u)\gamma_x$ is parallel to $\gamma$,
hence the presence of this term in (5.22) only amounts to a
reparametrization of the curve, and does not affect its shape.
The curve thus evolves in the direction of curvature.
An obvious Lyapunov functional is the {\it length} of the curve.
In terms of the variables
$$
\gamma_x=\pmatrix{v\cr w\cr}\doteq\pmatrix{u_x\cr -u_t\cr},\eqno(5.23)$$
this length is given by
$$L(\gamma)\doteq\int|\gamma_x|\,dx=\int\sqrt{v^2+w^2}\,dx\,.\eqno(5.24)$$
We can estimate the rate of decrease in the length as
$$-{d\over dt} L\big(\gamma(t)\big)~=~\int_{{I\!\!R}}
{|v|\,\big[(w/v)_x\big]^2\over \big(1+(w/v)^2\big)^{3/2}}\,dx
~\geq~
{1\over (1+\delta^2)^{3/2}}
\int_{|w/v|\leq \delta}
|v|\,\big[(w/v)_x\big]^2\,dx\,,\eqno(5.25)$$
for any given constant $\delta>0$.
This yields a useful a priori estimate on the
integral on the right hand side of (5.25).

\noindent{\bf 3.}
In connection with the same curve $\gamma$ in (5.21), we now introduce
another functional, defined in terms of a wedge product.
$$
Q(\gamma)\doteq{1\over 2}{\int\!\!\int}_{x<x'}\big| \gamma_x(x)\wedge
\gamma_x(x')\big|\,dx\,dx'\,.\eqno(5.26)$$
For any curve that moves in the plane in the direction of curvature,
one can show that
this functional is monotone decreasing and its decrease
bounds the area swept by the curve: $|dA|\leq -dQ$.

Using (5.22)--(5.23) we now compute
$$-{dQ\over dt}\geq \left|dA\over dt\right|=\int |\gamma_t\wedge
\gamma_x|\,dx=\int |\gamma_{xx}\wedge \gamma_x|\,dx=\int |v_x w-v
w_x|\,dx\,.$$ Integrating w.r.t.~time, we thus obtain another
useful a priori bound:
$$\int_0^\infty\int |v_x w-v w_x|\,dx\,dt~\leq~\int_0^\infty
\left| dQ\big(\gamma(t)\big)\over dt\right|\,dt~\leq~Q\big(\gamma(0)\big)\,.$$
Together, the functionals in (5.24) and (5.26)
allow us to
estimate all source terms in (5.12) due to the interaction of waves of the
same family.

This yields the ${\bf L}^1$ estimates on the source terms $\phi_i$,
in (5.12), proving the uniform bounds on the total variation of a solution
$u$ of (5.10).  See [BB] for details.
\vskip 0.8em
Next, to prove the uniform stability of all solutions
of the parabolic system (5.10) having small total variation, we
consider the linearized system describing the evolution of a first order
variation.
Inserting the formal expansion $u=u_0+\epsilon z +O(\epsilon^2)$
in (5.10), we obtain
$$z_t+\big[DA(u)\cdot z\big]u_x+A(u) z_x=z_{xx}\,.\eqno(5.27)$$
Our basic goal is to prove the bound
$$\big\|z(t)\big\|_{{\bf L}^1}\leq L\,
\big\|z(0)\big\|_{{\bf L}^1}\,,\eqno(5.28)$$
for some constant $L$ and all $t\geq 0$ and every solution $z$ of (5.27).
By a standard homotopy
argument, from (5.28) one easily deduces the Lipschitz continuity of the
solution of (5.8) on the initial data. Namely, for every couple of solutions
$u,\tilde u$ with small total variation one has
$$\big\| u(t)-\tilde u(t)\big\|_{{\bf L}^1}\leq L\,
\big\| u(0)-\tilde u(0)\big\|_{{\bf L}^1}\,.\eqno(5.29)$$
To prove (5.28)
we decompose the vector $z$ as a sum of scalar components:
$z=\sum_i h_i \tilde r_i$, write an evolution equation
for these components:
$$h_{i,t}+(\tilde\lambda_i h_i)_x-h_{i,xx}=\hat \phi_i\,,$$
and show that the source terms $\hat\phi_i$ are integrable on the domain
$\{ t>0\,, x\in{I\!\!R}\}$.
\vskip 0.8em
For every initial data $u(0,\cdot)=\bar u$ with small total variation,
the previous arguments yield the existence of a unique global solution
to the parabolic system (5.8), depending Lipschitz
continuously on the initial data, in the ${\bf L}^1$ norm.
Performing the rescaling $t\mapsto t/\varepsilon$, $x\mapsto x/\varepsilon$,
we immediately
obtain the same results for the Cauchy problem (5.2).
Adopting a semigroup notation, this solution can be written as
$u^\varepsilon(t,\cdot)=S^\varepsilon_t\bar u$.  Thanks to the uniform bounds on the
total variation, a compactness argument yields the existence of a strong
limit in ${\bf L}^1_{\rm loc}$
$$u=\lim_{\varepsilon_m\to 0}u^{\varepsilon_m}\eqno(5.30)$$
at least for some subsequence $\varepsilon_m\to 0$. Since the
$u^\varepsilon$ depend continuously on the initial data, with a
uniform Lipschitz constant, the same is true of the limit solution
$u(t,\cdot)=S_t\bar u$. In the conservative case where
$A(u)=Df(u)$, it is not difficult to show that this limit $u$
actually provides a weak solution to the Cauchy problem
(1.1)--(1.2).

The only remaining issue is to show that the limit in (5.30) is unique,
i.e.~it does not depend on the subsequence $\{\varepsilon_m\}$.  In the standard
conservative case, this fact can already be deduced from the
uniqueness result in [BG].
In the general case, uniqueness is proved in two steps.  First we
show that, in the special case of a Riemann problem,
the solution obtained as vanishing viscosity limit is unique
and can be completely characterized.
To conclude the proof, we then rely on
the same general argument as in [B2]: if two Lipschitz semigroups $S,S'$
provide the same solutions to all Riemann problems, then
they must coincide.  See [BB] for details.

\head{6. Concluding remarks}

\noindent{\bf 1.}
A classical tool in the analysis of first order hyperbolic systems
is the {\it method of characteristics}.
To study the system
$$u_t+A(u)u_x=0\,,$$
one decomposes the solution along the eigenspaces of the matrix $A(u)$.
The evolution of these components is then described by a family of
O.D.E's along the {\it characteristic curves}. In the $t$-$x$ plane,
these are the curves which satisfy $dx/dt=\lambda_i\big(u(t,x)\big)$.
The local decomposition (5.16)--(5.17) in terms of viscous travelling
waves makes it possible to implement this ``hyperbolic'' approach also
in connection with the parabolic system (5.10). In this case, the
projections are taken along the vectors $\tilde r_i$, while the
characteristic curves are defined as $dx/dt=\sigma_i$, where $\sigma_i$
is the speed of the $i$-th travelling wave. Notice that in the
hyperbolic case the projections and the wave speeds depend only on the
state $u$, through the eigenvectors $r_i(u)$ and the eigenvalues
$\lambda_i(u)$ of the matrix $A(u)$. On the other hand, in the parabolic
case the construction involves the derivatives $u_x$, $u_{xx}$ as well.
\vskip 0.8em \noindent{\bf 2.} In nearly all previous works on BV
solutions for systems of conservation laws, following [G] the basic
estimates on the total variation were obtained by a careful study of the
Riemann problem and of elementary wave interactions. The Riemann problem
also takes the center stage in all earlier proofs of the stability of
solutions [BC1], [BCP], [BLY]. In this connection, the hypothesis (H)
introduced by Lax [Lx1] is widely adopted in the literature. It
guarantees that solutions of the Riemann problem have a simple
structure, consisting of at most $n$ elementary waves (shocks, centered
rarefactions or contact discontinuities).  If the assumption (H) is
dropped, some results on global existence [L3], and continuous
dependence [AM] are still available, but their proofs become far more
technical. On the other hand, the approach introduced in [BB] marks the
first time where uniform $BV$ estimates are obtained without any
reference to Riemann problems. Global existence and stability of weak
solutions are obtained for the whole class of strictly hyperbolic
systems, regardless of the hypothesis (H). \vskip 0.8em \noindent{\bf
3.} For the viscous system of conservation laws
$$u_t+f(u)_x=u_{xx}\,,$$
previous results in [L4], [SX], [SZ], [Yu] have established the
stability of special types of solutions, for example travelling viscous
shocks or viscous rarefactions. Taking $\varepsilon=1$ in (5.2), from
Theorem~5 we obtain the uniform Lipschitz stability (w..r.t.~the ${\bf
L}^1$ distance) of ALL viscous solutions with sufficiently small total
variation. An interesting alternative technique for proving stability of
viscous solutions, based on spectral methods, was recently developed in
[HZ]. \vskip 0.8em \noindent{\bf 4.} In the present survey we only
considered initial data with small total variation.  This is a
convenient setting, adopted in much of the current literature, which
guarantees the global existence of $BV$ solutions of (1.1) and captures
the main features of the problem. A recent example constructed by
Jenssen [J] shows that, for initial data with large total variation, the
${\bf L}^\infty$ norm of the solution can blow up in finite time. In
this more general setting, one expects that the existence and uniqueness
of weak solutions, together with the convergence of vanishing viscosity
approximations, should hold locally in time as long as the total
variation remains bounded. For the hyperbolic system (1.1), results on
the local existence and stability of solutions with large $BV$ data can
be found in [Sc] and [BC2], respectively. Because of the counterexample
in [BS], on the other hand, similar well posedness results are not
expected in the general ${\bf L}^\infty$ case.

\parindent 12mm

\head{References}

\item{[AM]} F.~Ancona and A.~Marson, Well posedness for general
$2\times 2$ systems of conservation laws, {\it Amer. Math. Soc. Memoir},
to appear.

\item{[BaJ]} P.~Baiti and H.~K.~Jenssen, On the front tracking
algorithm, {\it J. Math. Anal. Appl.} {\bf 217} (1998), 395--404.

\item{[BB]}  S.~Bianchini and A.~Bressan, Vanishing viscosity solutions of
nonlinear hyperbolic systems,  preprint S.I.S.S.A., Trieste 2001.

\item{[B1]} A.~Bressan, Global solutions to systems of conservation laws by
wave-front tracking, {\it J. Math. Anal. Appl.} {\bf 170} (1992),
414--432.

\item{[B2]} A.~Bressan, The unique limit of the Glimm scheme,
{\it Arch. Rational Mech. Anal.} {\bf 130} (1995), 205--230.

\item{[B3]} A.~Bressan,
{\it Hyperbolic
Systems of Conservation Laws. The One Dimensional Cauchy Problem}.
Oxford University Press, 2000.

\item{[BC1]} A.~Bressan and R.~M.~Colombo, The semigroup generated
by $2\times 2$ conservation laws, {\it Arch. Rational Mech. Anal.}
{\bf 133} (1995), 1--75.

\item{[BC2]} A.~Bressan and R.~M.~Colombo, Unique solutions of $2\times 2$
conservation laws with large data, {\it Indiana Univ. Math. J.}
{\bf 44} (1995), 677--725.

\item{[BCP]} A.~Bressan, G.~Crasta and B.~Piccoli, Well posedness of the
Cauchy problem for $n\times n$ conservation laws, {\it
Amer. Math. Soc. Memoir} {\bf 694} (2000).

\item{[BG]} A.~Bressan and P.~Goatin, Oleinik type estimates and
uniqueness for $n\times n$ conservation laws, {\it J. Diff.
Equat.} {\bf 156} (1999), 26--49.

\item{[BLF]} A.~Bressan and P.~LeFloch, Uniqueness of weak solutions to
systems of conservation laws, {\it Arch. Rat. Mech. Anal.} {\bf
140} (1997), 301--317.

\item{[BLe]} A.~Bressan and M.~Lewicka, A uniqueness condition for
hyperbolic systems of conservation laws, {\it Discr. Cont. Dynam.
Syst.} {\bf 6} (2000), 673--682.

\item{[BLY]} A.~Bressan, T.~P.~Liu and T.~Yang,  $ L^1$ stability
estimates for $n\times n$ conservation laws, {\it Arch. Rational
Mech. Anal.} {\bf 149} (1999), 1--22.

\item {[BM]} A.~Bressan and A.~Marson, Error bounds for a
deterministic version of the Glimm scheme, {\it Arch. Rat. Mech.
Anal.} {\bf 142} (1998), 155--176.

\item{[BS]}
A.~Bressan and W.~Shen, Uniqueness for discontinuous O.D.E.~and
conservation laws, {\it Nonlinear Analysis, T. M. A.} {\bf  34}
(1998), 637--652.

\item{[C]} M.~Crandall, The semigroup approach to first-order quasilinear
equations in several space variables, {\it Israel J. Math.} {\bf
12} (1972), 108--132.

\item{[D1]} C.~Dafermos, Polygonal approximations of solutions of the initial
value problem for a conservation law, {\it J. Math. Anal. Appl.}
{\bf 38} (1972), 33--41.

\item{[D2]} C.~Dafermos, Hyperbolic systems of conservation laws,
{\it Proceedings of the International Congress of Mathematicians,
Z\"urich 1994}, Birch\'auser (1995), 1096--1107.

\item{[D3]} C.~Dafermos, {\it Hyperbolic Conservation Laws in Continuum
Physics}, Springer-Verlag, Berlin 2000.

\item{[DP1]} R.~DiPerna, Global existence of solutions to
nonlinear hyperbolic systems of conservation laws, { \it J. Diff.
Equat.}  {\bf 20} (1976), 187--212.

\item{[DP2]} R.~DiPerna, Convergence of approximate solutions to conservation
laws, {\it Arch. Rational Mech. Anal.} {\bf 82} (1983), 27--70.

\item{[G]} J.~Glimm, Solutions in the large for nonlinear hyperbolic systems
of equations, {\it Comm. Pure Appl. Math.} {\bf 18} (1965),
697--715.

\item{[GL]} J.~Glimm and P.~Lax, Decay of solutions of systems of nonlinear
hyperbolic conservation laws, {\it Amer. Math. Soc. Memoir} {\bf 101} (1970).

\item{[GX]} J.~Goodman and Z.~Xin, Viscous limits for piecewise smooth
solutions to systems of conservation laws, {\it Arch. Rational
Mech. Anal.} {\bf 121} (1992), 235--265.

\item{[HR]} H.~Holden and N.~H.~Risebro
{\it Front Tracking for Hyperbolic Conservation Laws},
Springer Verlag, New York 2002.

\item{[HZ]} P.~Howard and K.~Zumbrun, Pointwise semigroup methods
for stability of viscous shock waves, {\it Indiana Univ. Math. J.}
{\bf 47} (1998), 727--841.

\item{[K]} S.~Kruzhkov, First order quasilinear equations with several space
variables, {\it Math. USSR Sbornik} {\bf 10} (1970), 217--243.

\item{[J]}
H.~K.~Jenssen, Blowup for systems of conservation laws, {\it SIAM
J. Math. Anal.} {\bf 31} (2000), 894--908.

\item{[Lx1]} P.~Lax, Hyperbolic systems of conservation laws II, {\it
Comm. Pure Appl. Math.} {\bf 10} (1957), 537--566.

\item{[Lx2]} P.~Lax, Problems solved and unsolved concerning nonlinear P.D.E.,
{\it Proccedings of the International Congress of Mathematicians,
Warszawa 1983}. Elsevier Science Pub. (1984), 119--138.

\item{[L1]} T.~P.~Liu, The deterministic version of the Glimm scheme,
{\it Comm. Math. Phys.} {\bf 57} (1977), 135--148.

\item{[L2]} T.~P.~Liu, Linear and nonlinear
large time behavior of solutions of general systems of hyperbolic
conservation laws, {\it Comm. Pure Appl. Math.} {\bf 30} (1977),
767--796.

\item{[L3]} T.~P.~Liu, Admissible solutions of hyperbolic
conservation laws, {\it Amer. Math. Soc. Memoir} {\bf 240} (1981).

\item{[L4]} T.~P.~Liu, Nonlinear stability of shock waves,
{\it Amer. Math. Soc. Memoir} {\bf 328} (1986).

\item{[LY]} T.~P.~Liu and T.~Yang,
$L^1$ stability for $2\times 2$ systems of hyperbolic conservation
laws, {\it J. Amer. Math. Soc.} {\bf 12} (1999), 729--774.

\item{[M]} A.~Majda, {\it Compressible Fluid Flow and Systems of Conservation
Laws in Several Space Variables}, Springer-Verlag, New York, 1984.

\item{[O]} O.~Oleinik, Discontinuous solutions of
nonlinear differential equations (1957), {\it Amer. Math. Soc.
Translations} {\bf 26}, 95--172.

\item{[R]} N.~H.~Risebro, A front-tracking alternative to the random choice
method, {\it Proc. Amer. Math. Soc.} {\bf 117} (1993), 1125--1139.

\item{[Sc]} S.~Schochet, Sufficient conditions for local existence via
Glimm's scheme for large BV data, {\it J. Differential Equations}
{\bf 89} (1991), 317--354.

\item{[S1]} D.~Serre, {\it Systems of Conservation Laws I, II},
Cambridge University Press, 2000.

\item{[S2]} D.~Serre,
Systems of conservation laws : A challenge for the XXIst century,
{\it Mathematics Unlimited - 2001 and beyond},
B.~Engquist and W.~Schmid eds., Springer-Verlag, 2001.

\item{[SX]} A.~Szepessy and Z.~Xin, Nonlinear stability abd viscous shocks,
{\it Arch. Rational Mech. Anal.} {\bf 122} (1993), 53--103.

\item{[SZ]} A.~Szepessy and K.~Zumbrun, Stability of rarefaction waves in
viscous media, {\it Arch. Rational Mech. Anal.} {\bf 133} (1996),
249--298.

\item{[Yu]} S.~H.~Yu, Zero-dissipation limit of solutions with shocks for
systems of hyperbolic conservation laws, {\it Arch. Rational Mech.
Anal.} {\bf 146} (1999), 275--370.

\end